\documentclass[12pt,a4paper]{amsart}
\usepackage{amsmath,amssymb,amscd}
\input xy
\xyoption{all}
\newcommand{\lra}{\longrightarrow}

\newcommand{\ra}{\rightarrow}

\newcommand{\ZZ}{\mathbb Z}

\newcommand{\PP}{\mathbb P}

\theoremstyle{plain}
\newtheorem{theorem}{Theorem}[section]

\newtheorem{proposition}[theorem]{Proposition}

\newtheorem{remark}[theorem]{Remark}
\begin{document}
\title[Polarizations of Prym varieties]{Polarizations of Prym Varieties of pairs of coverings}

\author{H. Lange}
\author{S. Recillas}

\address{H. Lange\\Mathematisches Institut\\
              Universit\"at Erlangen-N\"urnberg\\
              Bismarckstra\ss e $1\frac{ 1}{2}$\\
              D-$91054$ Erlangen\\
              Germany}
\email{lange@mi.uni-erlangen.de}
\thanks{Supported by DAAD and Conacyt 40033-F}
\dedicatory{Sevin Recillas died on June 20, 2005.}

\begin{abstract}
To any pair of coverings $f_i: X \ra X_i, i = 1,2$ of smooth projective curves one can associate 
an abelian subvariety of the Jacobian $J_X$, the Prym variety $P(f_1,f_2)$ of the pair $(f_1,f_2)$. In some cases we 
can compute the type of the restriction of the canonical principal polarization of $JX$. We obtain 2 families of 
Prym-Tyurin varieties of exponent 6.

\end{abstract}

\maketitle

\section{Introduction}

Let $f: X \lra Y$ be a finite morphism of smooth projective curves. The complement of the abelian subvariety $f^*J_Y$ in the 
canonically polarized Jacobian $J_X$ of $X$ is called the Prym variety $P(f)$ of $f$. \noindent
In \cite{LR} we introduced an analogous notion for two morphisms of curves. To be more precise, 
suppose that we are given a commutative diagram of finite morphisms of smooth projective curves:
$$
\xymatrix@R=0.4cm@C=0.4cm{
 & X \ar[dl]_{f_1} \ar[dr]^{f_2} \\
         X_1 \ar[dr]_{g_1} && X_2 \ar[dl]^{g_2} \\
  & Y} \eqno(1.1)             
$$

\noindent
such that $g_1$ and $g_2$ do not both factorize via the same morphism $Y' \lra Y$  
of degree $\geq 2$. Then the {\it Prym variety $P(f_1,f_2)$ of the pair} $(f_1,f_2)$ is defined to be the complement
of the abelian variety $f_2^*(P(g_2))$ in $P(f_1)$ with respect to the canonical polarization. \\
By definition $P(f_1,f_2)$ is an abelian subvariety of the Jacobian
$J_X$ of the curve $X$.
In general it is difficult to determine the type of the restriction of a polarization to an abelian subvariety. 
It is the aim of this note to determine the type of the restiction of the canonical polarization of $J_X$ to $P(f_1,f_2)$ in some cases,
namely for $Y = \PP^1$ and $\deg f_1$ and $\deg f_2$ are prime to each other. As a special case we obtain 2 families of Prym-Tyurin varieties of exponent 6 
in any dimension $\geq 4$ respectively 5. Note that the correspondences defining $P(f_1,f_2)$ are not fixed point free. So one cannot apply Kanev's theorem
for the proof. \\    
In the last section we show that the Abel-Prym map of these abelian varieties is injective apart from the fact that 
all ramification points of $f_1$ are mapped to one point. Hence we obtain irreducible curves of low class in $P(f_1,f_2)$.
Moreover this implies that the Seshadri constant of $P(f_1,f_2)$ is small in 
these cases. 

\section{Restricting polarizations}

Let $(S,L)$ be a polarized abelian variety over an algebraically closed field of characteristic 0. As usual let $K(L)$ denote the
kernel of the associated homomorphism $\phi_{L}: S \ra \hat{S}, \,\, x \mapsto t_x^*L \otimes L^{-1}$. The polarization $L$ 
is of type $(d_1, \ldots, d_g)$ if and only if $K(L) \simeq (\ZZ/d_1\ZZ \times \cdots \times \ZZ/d_g\ZZ)^2$.\\

Now let $A$ be an abelian subvariety of $S$ with canonical embedding $\iota_A : A \lra S$. 
Let 
$$
P := \ker(\hat{\iota}_A \circ \phi_L: S \lra \hat{A})^0   \eqno(2.1)
$$
denote the complementary abelian subvariety of $S$. Let $L_A := L|A$ and $L_P := L|P$ denote the restriction of the polarization $L$ 
to $A$ and $P$. If $L$ is a principal polarization, it is well known that 
$$
K(L_A) \simeq K(L_P) \simeq A \cap P                \eqno(2.2)
$$ 
from which the type of $L_A$ and $L_P$ can be deduced in many cases. For arbitrary polarizations the situation is more complicated.
In this section we recall two generalizations of (2.2) which will be applied later and for the proof of which we refer to 
\cite{LP}.\\ 
 
\begin{proposition}
$| K(L_P)| \cdot |K(L_A)| = | A \cap P|^2 \cdot |K(L)|.$
\end{proposition} 

\begin{proposition}
There is an exact sequence 
$$
0 \lra K(L) \cap P \lra K(L_P) \lra A \cap P \lra 0.
$$
\end{proposition}

\noindent
In particular $|K(L_P)| = |K(L) \cap P| \cdot |A \cap P|$ and similarly for $L_A$.\\

\section{Restriction of the canonical polarization to $P(f_1,f_2)$}

Let $f: X \lra Y$ be a morphism of degree $n$ of smooth projective curves over an algebraically closed field $k$.
Denote by $J_X := \mathrm{Pic} ^0 (X)$ and $J_Y := \mathrm{Pic} ^0 (Y)$ the Jacobians of $X$ and $Y$. Pulling back line bundles 
defines a homomorphism 
$$
f^* : J_Y \lra J_X.
$$
$f^*$ has finite kernel and is an embedding if and only if $f$ does not factor via a cyclic \'etale cover 
of degree $\geq 2$ (see \cite{CAV}, Proposition 11.4.3). The norm map of line bundles defines 
a homomorphism
$$
N_f
: J_X \lra J_Y.
$$
The {\it Prym variety $P(f)$ of the morphism} $f$ is defined to be the abelian subvariety
$$
P(f) := \ker (N_f)^0
$$
of $J_X$ where the $0$ means the connected component containing $0$.\\ 

\noindent
Now suppose that we are given a commutative diagram of finite morphisms of smooth projective curves:
$$
\xymatrix@R=0.4cm@C=0.4cm{
 & X \ar[dl]_{f_1} \ar[dr]^{f_2} \\
         X_1 \ar[dr]_{g_1} && X_2 \ar[dl]^{g_2} \\
  & \PP^1             } 
  \eqno(3.1)
$$
with the following properties:
\begin{itemize}
\item[(a)] $g_1$ and $g_2$ do not both factorize via the same morphism $Y' \ra \PP^1$  
 of degree $\geq 2$,
\item[(b)] $X$ is the normalisation of the cartesian product $X_1 \times_{\PP^1} X_2$. 
\item[(c)] $d_1 := \deg g_1$ and $d_2 := \deg g_2$ are prime to each other. 
\end{itemize}

\noindent
Then we have for the Prym variety $P = P(f_1,f_2)$, defined as the complement of $f_2^*(P(g_2))$ in $P(f_1)$ 
with respect to the restriction of the canonical polarization of $J_X$ to $P(f_1)$ (see \cite{LR}):\\

\begin{proposition}
Let $\Theta$ be the canonical principal polarization of $J_X$ and denote
$$
L := \Theta|P(f_1), \,\, L_{f_2^*P(g_2)} := L|f_2^*P(g_2) = \Theta|f_2^*P(g_2), \,\, L_P := L|P = \Theta|P.
$$
Then
$$
K(L_P) \simeq K(L) \oplus K(L_{f_2^*P(g_2)}).
$$
\end{proposition}

\begin{proof}
$P(f_1)$ is by definition the complementary abelian subvariety of $f_1^*J_{X_1}$ in $J_X$ with repect to 
$\Theta$. According to \cite{CAV}, Proposition 11.4.3 $\Theta|f_1^*J_{X_1}$ is of type $(d^{(1)}_2, \cdots ,d^{(s)}_2,d_2, \cdots, d_2)$ where 
$d^{(i)}_2$ divides $d_2$ for $i = 1, \ldots, s$ and $0 \leq s < \dim J_{X_1}$. Hence, according to (2.2) the polarization $L$ is of type
$(1, \ldots,1,d^{(1)}_2, \cdots ,d^{(s)}_2,d_2, \cdots, d_2)$. In particular $K(L)$ is of exponent $d_2$.\\
On the other hand $P(g_2) = J_{X_2}$. Again by \cite{CAV}, Proposition 11.4.3, $L_{f_2^*P(g_2)} = \Theta|f_2^*P(g_2)$ is of type 
$(d^{(1)}_1, \cdots ,d^{(t)}_1,d_1, \cdots, d_1)$ where $d^{(i)}_1$ divides $d_1$ for $i = 1, \ldots, t$ and $0 \leq t < \dim J_{X_2}$. 
In particular $L_{f_2^*P(g_2)}$ is of exponent $d_1$.\\
Now $P$ is by definition the complementary abelian subvariety of $f_2^*P(g_2)$ in $P(f_1)$ with respect to the polarization $L$. 
Hence according to Proposition 2.2 there is an exact sequence
$$
0 \lra K(L) \cap f_2^*P(g_2) \lra K(L_{f_2^*P(g_2)}) \lra  P \cap f_2^*P(g_2)  \lra 0
$$
The assumption $gcd(d_1,d_2) = 1$ implies by the above reasoning that $gcd(|K(L)|,|K(L_{f_2^*P(g_2)}|) = 1$.
Hence $K(L) \cap f_2^*P(g_2) = 0$ and thus 
$$
K(L_{f_2^*P(g_2)}) \simeq f_2^*P(g_2) \cap P.
$$
But then Proposition 2.1 yields 
$$
|K(L_P)| = |K(L)| \cdot |K(L_{f_2^*P(g_2)})|
$$
which implies the assertion, since $d_1$ and $d_2$ are prime to each other.
\end{proof}

\noindent
By symmetry reasons we also get with the same notation\\

\begin{remark} 
{\it Denote
$$
M:= \Theta|P(f_2), \quad  M_{f_1^*P(g_1)} := M|f_1^*P(g_1), \quad M_P = \Theta|P = L_P.
$$
Then}
$$
K(L_P) = K(M_P) \simeq K(M) \oplus K(M_{f_1^*P(g_1)}).
$$
\end{remark}

Recall that $P$ is a Prym-Tyurin variety of exponent $q$ in the Jacobian $J_X$ if the restricted polarization $L_P = \Theta|P$ is 
the $q$-fold of a principal polarization of $P$. This is the case if and only if $L_P$ is of type $(q, \ldots, q)$.

\begin{proposition}
Let the notation be as above. In particular we assume $gcd(d_1,d_2) = 1$. If we assume moreover that $f_1$ and $f_2$ 
do not factorize via a cyclic unramified covering, then the following conditions are equivalent
\begin{itemize}
\item[(a)] $P = P(f_1,f_2)$ is a Prym-Tyurin variety of exponent $d_1d_2$ in $J_X$.
\item[(b)] $\dim(P) = g(X_1) = g(X_2)$.
\end{itemize}
\end{proposition} 

\begin{proof}
Since $f_1$ does not factorize via an unramified covering $L$ is of type $(1,\ldots,1, d_2, \ldots, d_2)$ with $g(X_1)$ entries $d_2$.
Since $f_2$ does not factorize via an unramified covering, $L_{f_2^*P(g_2)}$ is of type $(d_1, \ldots, d_1)$ with $g(X_2)$ entries $d_1$.
Hence according to Proposition 3.1 $L_P$ is of type $(1,\ldots,1, d_1d_2, \ldots, d_1d_2)$ if $g(X_1) = g(X_2)$ with $g(X_1)$ enties $d_1d_2$ and 
$\dim P - g(X_1)$ entries 1. In particular we have $\dim P \geq g(X_1)$. If $\dim P = g(X_1) = g(X_2)$, there is no entry 1 implying that 
$(P,L_P)$ is a Prym-Tyurin variety of exponent $d_1d_2$ in $J_X$. Conversely, if one of the equations in (b) is not satisfied, the above reasoning 
implies that $(P,L_P)$ cannot be a Prym-Tyurin variety in $J_X$. 
\end{proof}

We also need the following variation of Proposition 3.3:

\begin{proposition}
Let the notation be as in Propossition 3.3 and assume now that as above $f_1$
does not factorize via an unramified covering, but $f_2$ is \'etale cyclic. Then the following conditions are equivalent
\begin{itemize}
\item[(a)] $P = P(f_1,f_2)$ is a Prym-Tyurin variety of exponent $d_1d_2$ in $J_X$.
\item[(b)] $\dim(P) = g(X_1) = g(X_2) - 1 \,\,\, \mbox{and} \,\,\, d_1 = 3$.
\end{itemize}
Moreover in this case $d_1=3$ if $\dim P \geq 2$.
\end{proposition} 

\begin{proof}
As above $L$ is of type $(1,\ldots, 1, d_2, \ldots, d_2)$ with $g(X_1)$ entries $d_2$. Since $f_2$ is \'etale cyclic of degree $d_1$, 
$L_{f_2^*P(g_2)}$ is of type $(1, d_1, \ldots, d_1)$ with $g(X_2) - 1$ entries $d_1$.
Hence according to Proposition 3.1 $L_P$ is of type $(1,\ldots,1, d_1d_2, \ldots, d_1d_2)$ if $g(X_1) = g(X_2) - 1$ with $g(X_1)$ entries $d_1d_2$ and 
$\dim P - g(X_1)$ entries 1. As in the proof of Proposition 3.3 one concludes that if $\dim(P) = g(X_1) = g(X_2) - 1$, then 
$P$ is a Prym-Tyurin variety.\\
Moreover, since $f_2$ is cyclic \'etale, we have $g(X) = d_1 g(X_2) - d_1 + 1$ and thus
\begin{eqnarray*}
g(X_2) - 1 = \dim P & = & g(X) - g(X_1) - g(X_2)\\
& = & (d_1 - 2) g(X_2) - d_1 + 2, 
\end{eqnarray*}
implying $d_1 = 3$ if $g(X_2) = \dim P + 1 \geq 3$. Again, if the equation in (b) are not satisfied, it is easy to see that $(P,L_P)$ cannot be a Prym-Tyurin variety in $J_X$.
\end{proof}

\section{ Two families of Prym-Tyurin varieties}

In this section we analyse the conditions of Propositions 3.3 and 3.4 to obtain two families of Prym varieties of pairs of morphisms which are Prym-Tyurin
varieties. Suppose we are given a commutative diagram (3.1). For any finite morphism $f$ of smooth projective curves let $\delta(f)$ denote the degree of 
the ramification divisor. This is always an even number. Hence we may denote
$$
\delta(g_i) = 2 r_i \qquad \mbox{and} \qquad \delta (f_i) = 2 s_i \qquad \mbox{for} \qquad i = 1, 2.
$$
Of course $r_1, r_2, s_1$ and $s_2$ are not unrelated. Then Hurwitz' formula gives
$$
g(X_i) = r_i - d_i + 1  \,\,\, \mbox{for} \,\,\, i = 1,2 
$$  
and 
$$
g(X) = d_2(r_1 - d_1) + s_1 + 1 = d_1(r_2 - d_2) + s_2 + 1.
$$
This implies
\begin{eqnarray*}
\dim P & = & d_1r_2 - d_1d_2 +d_1 +d_2 -r_1 -r_2 +s_2 -1\\
& = & d_2r_1 - d_1 d_2 +d_1 + d_2 - r_1 - r_2 + s_1 - 1
\end{eqnarray*}
In particular we have
$$
s_2 = d_2r_1 - d_1r_2 + s_1.   \eqno(4.1)
$$

\begin{theorem}
Let the assumptions be as in Proposition 3.3. In particular we assume that $gcd(d_1,d_2) = 1$ and that $f_1$ and $f_2$ do not factorize via a
cyclic \'etale covering. Without loss of generality we may assume $d_1 > d_2$. \\
Then $P = P(f_1,f_2)$ is a Prym-Tyurin variety of dimension $ \geq 4$ in $J_X$ if and only if
$$
d_1 = 3, d_2 = 2, r_1 \geq 6, s_1 = r_2 = r_1 - 1 \,\,\, \mbox{and} \,\,\, s_2 = 2. 
$$ 
In this case $P$ is a Prym-Tyurin variety of dimension $r_1 - 2$ in $J_X$.
\end{theorem}

\begin{proof}
According to Proposition 3.3 $P$ is a Prym-Tyurin variety of exponent $d_1d_2$ in $J_X$ if and only if $\dim P = g(X_1) = g(X_2)$ which is equvalent to\\
$$
\left\{ \begin{array}{ccc} 
r_1 - d_1 & = & r_2 - d_2\\
d_1r_2 - d_1d_2 +d_1 +d_2 -r_1 -r_2 +s_2 -1 & = & r_2 - d_2 + 1
\end{array} \right.
$$
Using also (4.1) this is the case if and only if the following 3 equations hold
$$
r_2  =  r_1 + d_2 - d_1    \eqno(4.2)
$$
$$
s_1  = (3-d_2)r_1 + (d_2 - 3)d_1 + 2     \eqno(4.3)
$$
$$
s_2  =  (3-d_1)r_1 + (d_1 - 3)d_1 + 2     \eqno(4.4)
$$\\

{\it Case} 1: $d_1 > d_2 \geq 4$.\\

By (4.4) and since $d_1 \geq 5$ we have $r_1 \leq d_1 + \frac{2}{d_1 -3} \leq d_1 + 1$. On the other hand (4.2) 
implies $r_1 \geq d_1 - d_2 + 1$.
Hence we obtain
$$
r_1 = d_1 - d_2 + i   \,\,\, \mbox{with} \,\,\, 1 \leq i \leq d_2 + 1
$$
and thus 
\begin{eqnarray*} 
\dim P & = & d_2r_1 - d_1d_2 + d_1 + d_2 - r_1 - r_2 + s_1 - 1\\
& = & i - d_2 + 1 \leq 2
\end{eqnarray*}
Hence there is no Prym-Tyurin variety $P$ of dimension $\geq 3$ in this case.\\

{\it Case} 2: $d_1 \geq 4,  \,\, d_2 = 3$.\\

Here we have by (4.4) $r_1 \leq d_1 + \frac{2}{d_1 - 3} \leq d_1 + 2$, whereas (4.2) gives $r_1 \geq d_1 - 3$. Hence we get
$$
r_1 = d_1 - 3 + i \,\,\, \mbox{with} \,\,\, 0 \leq i \leq 5
$$
and thus using that $s_1 = 2$ by (4.3),
\begin{eqnarray*}
\dim P & = & d_2r_1 - d_1d_2 + d_1 + d_2 - r_1 - r_2 + 1\\
& = & r_1 - d_1 + 1 = i - 2 \leq 3.
\end{eqnarray*}
Hence in this case there is no Prym-Tyurin variety $P$ of dimension $\geq 4$.\\

{\it Case} 3: $d_1 \geq 4, \,\, d_2 = 2$.\\

Again we have by (4.4) $r_1 \leq d_1 + \frac{2}{d_1 - 3}$. But now $d_1 \geq 5$ since $(d_1,d_2) = 1$, implying $r_1 \leq d_1 + 1$.
As above we get $r_1 \geq d_1 - 3$ and hence
$$
r_1 = d_1 -3 + i \,\,\, \mbox{with} \,\,\, 0 \leq i \leq 4
$$  
and thus using that $s_1 = r_2$ by (4.2) and (4.3)
\begin{eqnarray*}
\dim P & = & r_1 - d_1 + 1\\
& = & i - 2 \leq 2.
\end{eqnarray*} 
Hence in this case there no Prym-Tyurin variety of dimension $\geq 3$.\\

{\it Case} 4: $d_1 = 3, \,\, d_2 = 2$.\\

Here equations (4.2) - (4.4) say $s_1 = r_2 = r_1 - 1$ and $s_2 = 2$ and thus are fulfilled for any $r_1 \geq 2$. According to what 
we have said above $P$ is a Prym-Tyurin variety of exponent 6 in $J_X$. Its dimension is $\dim P = r_1 - 2$.  
\end{proof}

\begin{remark} In Theorem 4.1 we obtain for every dimension $d \geq 4$ a family of Prym-Tyurin varieties of dimension d and exponent 6. Counting parameters 
it is easy to see that the family is of dimension $2d + 1$.
\end{remark}

\begin{theorem}
Let the assumptions be as in Proposition 3.4. In particular we assume that $gcd(d_1,d_2) = 1$, that $f_1$ does not factorize via a
cyclic \'etale covering and that $f_2$ is \'etale cyclic. \\
Then $P = P(f_1,f_2)$ is a Prym-Tyurin variety of dimension $ \geq 5$ in $J_X$ if and only if
$$
d_1 = 3, \,\, d_2 = 2, \,\, r_1 = r_2 = s_1 \geq 7 \,\, \mbox{and} \,\,   s_2 = 0.
$$
In this case $P$ is a Prym-Tyurin variety of dimension $r_1 - 2$ and of exponent 6 in $J_X$.
\end{theorem}

\begin{proof}
According to Proposition 3.4 $P$ is a Prym-Tyurin variety of exponent $d_1d_2$ in $J_X$ if and only if $\dim P = g(X_1) = g(X_2) - 1$ and  
$d_1 = 3$. Using (4.1) this is the case if and only if
$$
\left\{ \begin{array}{c}
d_1 = 3 \,\,\, \mbox{and} \,\,\,  s_2 = 0\\
r_2 = d_2 + r_1 - 2\\
s_1  = 3 r_2 - d_2 r_1 
\end{array} \right. \eqno(4.5)
$$   
Now $s_1 \geq 0$ and hence (4.5) implies
$$
(d_2 - 3)(3 - r_1) + 3 = 3 d_2 + 3 r_1 -6 - d_2 r_1  \geq 0
$$
On the other hand $r_1 > 0$ and $r_2> 0$, since there are no nontrivial unramified coverings of $\PP^1$.\\
If $d_2 \geq 7$, this gives $r_1 \leq 3$ and thus $\dim P = r_1 - 2 \leq 1$. If $d_2 = 6$ or $5$, then $r_1 \leq 4$ 
and thus $\dim P = r_1 - 2 \leq 2$. If $d_2 = 4$, then $r_1 \leq 6$ and thus $\dim P = r_1 -2 \leq 4$. Hence 
in these cases there is no Prym-Tyurin variety $P$ of dimension $\geq 5$.\\
We are left with $d_2 = 2$ which leads to the conditions of the theorem. According to Proposition 3.4 
in this case $P$ is a Prym-Tyurin variety of dimension $r_1 - 2$.  
\end{proof}

\begin{remark}
It is easy to see that in the case of Theorem 4.3 the composed covering $g_i \circ f_i: X \ra \PP^1$ is Galois with Galois group $S_3$.
Again we obtain for every dimension $d \geq 5$ a family of Prym-Tyurin varieties of dimension d and exponent 6. Counting parameters
it is easy to see that the family is of dimension $2d + 1$.
\end{remark}

\begin{remark}
It is easy to see that these are the only cases of pairs of morphisms $(f_1, f_2)$ as in diagram (3.1) where we get a 
Prym-Tyurin variety in this way.
\end{remark}

\section{ The Abel-Prym map}

\noindent
Let the notation be as in Theorems 4.1 or 4.3. In particular $P = P(f_1,f_2)$ is the Prym-Tyurin variety of the pair $(f_1,f_2)$.
Let $\alpha: X \ra J_X$ denote the Abel map with respect to a base point $x_0$. Identifying $J_X$ and $P$ with their dual abelian 
varieties via the principal polarizations, the canonical embedding $\iota_P:P \ra J_X$ induces a surjective homomorphism 
$\hat{\iota}: J_X \ra P$. The composed map
$$
\beta: X \stackrel{\alpha}{\ra} J_X \stackrel{\hat{\iota}}{\ra} P
$$
is called the {\it Abel-Prym map} of the Prym-Tyurin variety $P$. Let $R_1 \subset X$ denote the branch locus of the double covering
$f_1:X \ra X_1$.  

\begin{proposition}
If $\dim P \geq 8$, then we have for the Abel-Prym map $\beta: X \ra P$:\\ $\beta$ is injective
on $X \setminus R_1$ and contracts $R_1$ to a point in $P \setminus \beta(X \setminus R_1)$. 
\end{proposition}

\begin{proof}
The Prym-Tyurin variety $P$ is the complement of the abelian subvariety $f_1^*J_{X_1} + f_2^*J_{X_2}$ of $J_X$ with respect to the 
canonical polarization. $f_1^*J_{X_1} + f_2^*J_{X_2}$ is of exponent 6 in $J_X$. Let $N_{f_i}: J_X \ra J_{X_i}$ denote the Norm map for $i=1,2$. 
It is easy to see that
$$
\hat{\iota} = 6\cdot 1_{J_X} - 3N_{f_1} - 2N_{f_2}
$$
Assume that there are points $x \not= y$ of $X$ such that $\beta(x) = \beta(y)$. If $\iota_1$ 
denotes the involution of $X$ induced by $f_1$ and for any $x \in X$ we denote by $\iota_2(x)$ and $\iota_2^2(x)$
the other 2 points of the fibre $f_2^{-1}(f_2(x))$ (note that $\iota_2$ is not a global automorphism of $X$, but locally this
notation makes sense), this means 
$$
6x -3(1 +\iota_1)(x) -2(1 + \iota_2 - \iota_2^2)(x) \sim 6y -3(1 +\iota_1)(y) -2(1 + \iota_2 - \iota_2^2)(y)
$$
or equivalently
$$
x + 3\iota_1(y) + 2\iota_2(y) + 2\iota_2^2(y) \sim y + 3\iota_1(x) + 2\iota_2(x) + 2\iota_2^2(x) \eqno(6.1)
$$
Here $\sim$ means linear equivalence of divisors on the curve $X$.\\
It suffices to show that if $x$ or $y \notin R_1$ this implies that (6.1) defines a pencil of divisors of degree 8, since then \\
either the $g_8^1$ does not factor via $f_1$ in which case Casteluovo's formula gives $g(X) \leq 1\cdot 7 + 2g(X_1)$ giving
$3r_1 - 6 \leq 7 + 2r_1 - 4$ which implies $r_1 \leq 9$,\\
or the $g_8^1$ factors via $f_1$ in which case Castelnuovo's formula applied to $X_1$ gives $r_1 - 2 = g(X_1) \leq 2\cdot 3$ which implies $r_1 \leq 8$.
But by assumption $ r_1 = \dim P + 2\geq 10$, a contradiction.\\
In order to complete the proof, assume (6.1) is an equality of divisors. 

Assume first that $x$ and $y$ are in the same fibre of 
the composed map $g_i\circ f_i: X \ra \PP^1$. Since certainly $3\iota_1(y) \not= 3\iota_1(x)$ we may assume after eventually exchanging
$\iota_2(x)$ and $\iota_2^2(x)$ respectively $\iota_2(y)$ and $\iota_2^2(y)$ that
$$
3\iota_1(y) = y + 2\iota_2(x) \quad \mbox{and} \quad 3\iota_1(x) = x + 2\iota_2(y) 
$$
In particular $x$ and $y$ are two different ramification points of $f_1$ in a fibre of $g_1 \circ f_1$. But every fibre of  
$g_1 \circ f_1$ contains at most one ramification point of $f_1$, a contradiction.\\
If $x$ and $y$ are not in the same fibre of $g_i\circ f_i$ and $f_1$ is not ramified at $x$ or $y$, then (6.1) cannot 
be an equality of divisors and the same arguement works. So finally assume that $x$ and $y$ are different points of $R_1$. 
Then (6.1) reads 
$$
2y + 2\iota_2(y) + 2\iota_2^2(y) \sim 2x + 2\iota_2(x) + 2\iota_2^2(x)
$$ 
which means that (6.1) is the $g_6^1$ defined by $g_i\circ f_i$. In particular $\beta(R_1)$ contracts to a point. 
\end{proof}

\begin{remark}
Observe that with the notation of the proof of Proposition 5.1 the correspondence $D_1$ of the curve $C$ defined by 
$$
D_1(x) = x -3\iota_1(x) - 2 \iota _2(x) - 2\iota_2^2(x)
$$
for any $x \in X$ induces the endomorphism $\hat{\iota}$ with image $P$. The correspondence $D$ of $C$ defined by
$$
D(x) = \pi^*\pi(x) - D_1(x)  
$$
is an effective correspondence defining $\hat{\iota}$. But $D$ is not fixed point free, so one cannot apply Kanev's 
theorem (see \cite{CAV}, Theorem 12.9.1) in order to show that $P$ is a Prym-Tyurin variety. 
\end{remark}

\begin{remark}
One can show that the singularity $\beta(R_1)$ of the curve $\beta(X)$ is an ordinary singularity. 
We do not include the details here.
\end{remark}

\begin{remark}
Recall that the Seshadri constant of the polarisation $\Xi$ of the abelian variety $P$ is defined as
$$
\varepsilon(\Xi) = \inf \frac{(\Xi \cdot C)}{mult_0(C)}
$$
where the multiplicity is to be taken over all irreducible reduced curves $C$ on $P$ containing 0. The 
existence of a curve with high multiplicity means that the Seshadri constant is small.
\end{remark}

\noindent
In fact 
$$
\varepsilon(\Xi) \leq \left\{ \begin{array}{c} 
3 - \frac{5}{\dim P +1}\\
3 - \frac{5}{\dim P + 2}
\end{array} \right.
\mbox{if} \quad P \quad \mbox{is as in} 
\left\{ \begin{array}{c} 
\mbox{Theorem} \quad 4.1,\\
\mbox{Theorem} \quad 4.3.
\end{array} \right.
$$
This follows from the fact that $C$ is of class $\frac{6}{(\dim P -1)!}[\Xi]^{\dim P -1}$.

\end{document}